\documentclass[10pt,leqno,twoside]{article}

\usepackage{annal-latex}
\usepackage{graphicx}

\usepackage{url}
\usepackage[pdftex, colorlinks, citecolor=blue, linkcolor=red, urlcolor=blue]{hyperref}%
\usepackage{enumerate}

\newtheorem{theorem}{\noindent Theorem}[section]

\newtheorem{definition}{\noindent Definition}[section]
\newtheorem{example}[theorem]{\noindent Example}
\newtheorem{remark}{\noindent Remark}[section]
\newtheorem{conjecture}{\noindent Conjecture}

\newcommand*{\Var}{\mathop{\mathrm{Var}}}

\newcommand*{\dCov}{\mathop{\mathrm{dCov}}}
\newcommand*{\dCor}{\mathop{\mathrm{dCor}}}
\newcommand*{\eCov}{\mathop{\mathrm{eCov}}}
\newcommand*{\eCor}{\mathop{\mathrm {eCor}}}
\newcommand*{\eVar}{\mathop{\mathrm{eVar}}}
\newcommand*{\E}{\mathbb E}
\renewcommand*{\Pr}{\mathbb P}

\begin{document}
\setcounter{page}{1}
\newcommand\balline{\small Tam\'as F.\ M\'ori and G\'abor J.\ Sz\'ekely}
\newcommand\jobbline{\small The Earth Mover's Correlation}

\vspace{-4cm} \fofej{50}{20}{nn}{nnn}

\vspace{.4cm}

\title{THE EARTH MOVER'S CORRELATION}

\author{{\bf Tam\'as F. M\'ori} (Budapest, Hungary)\\[1ex]
{\bf G\'abor J. Sz\'ekely}(Alexandria, VA, USA)
}

%
\mathclass{62H20}

\projsupport{T.\  F.\ M\'ori was supported by the Hungarian National 
Research, Development and Innovation Office NKFIH -- Grant No. K125569. 
Part of this research was based on work supported by the National 
Science Foundation, while the second author was working at the
Foundation. G.\ J.\ Sz\'ekely is grateful for many interesting discussions 
with G\'abor Tusn\'ady, Endre Boros, Imre Z. Ruzsa, and Russell Lyons.}
%



\commby{???}
\vspace{-2ex}
\recacc{???}{???}

\vspace{-7ex}

\abstract {Since Pearson's correlation was introduced at the end of the 19th century many 
dependence measures have appeared in the literature. In  \cite{ms2019} we suggested 
four simple axioms for dependence measures of random variables that take values in 
Hilbert spaces. We showed that distance correlation (see \cite{srb07}) satisfies all these axioms. 
We still need a new measure of dependence because existing measures either do not 
work in general metric spaces (that are not Hilbert spaces) or they do not satisfy our four 
simple axioms. The {\it earth mover's correlation} introduced in this paper applies in general 
metric spaces and satisfies our four axioms (two of them in a weaker form).}

\section{Introduction: What is our goal?}
\label{intro}

Let $S$ be a nonempty set of pairs of nondegenerate random variables  $X,Y$ 
taking values in Euclidean spaces or in real, separable Hilbert spaces $H$. 
(Nondegenerate means that the random variable is not constant with probability $1$.) 
In \cite{ms2019}  we called $\Delta (X, Y): S \to [0,1]$  a dependence measure on $S$
if the following four axioms hold. 

In the axioms below we need similarity transformations of $H$. Similarity 
$H$ is defined as a bijection (\textup{1--1} correspondence) from $H$ 
onto itself that multiplies all distances by the same positive real number (scale). 
Similarities in Hilbert spaces are known to be compositions of a translation, an 
orthogonal linear mapping, and a uniform scaling. We assume that if $(X, Y) \in S$ 
then $\big(f(X),g(Y)\big) \in S$ for all similarity transformations $f,g$ of $H$. 

In \cite{ms2019} we introduced the following axioms. 
\begin{enumerate}[(i)]
\item\label{ni} $\Delta (X, Y) = 0$ if and only if $X$ and $Y$ are
  independent. 
\item\label{nii} $\Delta (X, Y)$ is invariant with respect to all
  similarity transformations of $H$; that is, $\Delta\big(f(X), g(Y)\big) =
  \Delta(X,Y)$ where $f$, $g$ are similarity transformations of $H$.
\item\label{niii} $\Delta (X, Y) = 1$ if and only if $Y = f(X)$  with
  probability $1$, where $f$ is a similarity transformation of $H$.
\item\label{niv} $\Delta (X, Y)$ is continuous; that is, if for some
  positive constants $K$ we have $E\big(|X_n|^2+|Y_n|^2\big)\le K$, $n=1,2,\dots$
  and $(X_n, Y_n)$ converges weakly (converges in distribution) to $(X,Y)$
  then $\Delta (X_n, Y_n) \to \Delta(X,Y)$.
\end{enumerate}

In fact, what we really need is not the boundedness of the second moments but the 
convergence of the expectations: $\E(X_n) \to \E(X)$ and $\E(Y_n) \to \E(Y)$;
so in axiom (\ref{niv}) the condition on the boundedness of second moments can be 
replaced by any other condition that guarantees the convergence of expectations. 
Such a condition is uniform integrability of $X_n, Y_n$ which follows from the 
boundedness of second moments. The reason of using a more restrictive condition is 
that it can be more easily checked.

If $S$ is the set of bivariate Gaussian random variables then Pearson's correlation 
satisfies all these axioms. For more general $S$ Pearson's correlation typically does 
not satisfy (\ref{ni}) but distance correlation does satisfy all of them if the expectations are finite.

First of all recall the definition of the sample distance correlation, see also \cite{srb07} 
and \cite{sr09}. Take all pairwise distances between sample values of one variable, and
do the same for the second variable. Rigid motion invariance is 
automatically guaranteed if instead of sample elements we work with 
their distances. Another advantage of working with distances is that 
they are always real numbers even when the data are vectors of 
possibly different dimensions. Once we have computed the distance 
matrices of both samples, double-center them (so each has column and 
row means equal to zero). Then average the entries of the matrix which 
holds componentwise products of the two centered distance 
matrices. This is the square of the sample distance covariance. 
If we denote the centered distances by $A_{ij}$, $i,j =1,\dots, n$ and
$B_{ij}$, $i,j =1,\dots, n$ where $n$ is the sample size, then the 
squared sample distance covariance is 
\[
\frac {1}{n^2} \sum_{i,j=1}^n A_{i,j}B_{i,j}.
\]
This definition is very similar to, and almost equally simple as, the
 definition of Pearson's covariance, except that here we have double 
indices.

If $\E|X|^2$ and $\E|Y|^2$ are finite then the population squared distance 
covariance can be reduced to the following form \cite{srb07}. 

If $(X,Y)$, $(X',Y')$, $(X'',Y'')$ denote independent and identically distributed copies then
 the distance covariance is the square root of
\begin{multline*} 
{\dCov}^2(X,Y):= \E(|X -X'|\,|Y -Y'|) + \E(|X-X'|)\E(|Y-Y'|)\\  
- \E(|X-X'|\,|Y-Y''|) -\E(|X -X''|\,|Y -Y'|).
\end{multline*}

In the above referenced paper we proved that the distance variance, 
$\dCov(X,X)$ is zero if and only if $X$ is constant with probability $1$. 
Once we defined distance covariance and distance variance we can define distance 
correlation the same way as we defined correlation with the help of 
covariance and variance. If the random variables $X,Y$ have finite 
expected values and they are not constant with probability $1$ then 
the definition of {\it population distance correlation} is the following:
\[
\dCor(X,Y):=\frac{\dCov(X,Y)}{\sqrt{\strut\dCov(X,X)\dCov(Y,Y)}} . 
\]
If $\dCov(X,X)\dCov(Y,Y)= 0$ then we do not define $\dCor(X,Y)$. 
(If we define $\dCor(X,Y) = 0$ then this would lead to a violation of (\ref{niv}).)

Distance correlation equals zero if and only if the 
variables are independent, whatever be the underlying distributions 
and whatever be the dimension of the two variables (for a transparent 
explanation see below). This fact and the 
simplicity of the statistic make distance correlation an
 attractive candidate for measuring dependence. 
For generalizations to certain metric spaces see \cite{lyons2013}, \cite{lyons2018}, and \cite{Jakobsen}.
These metric spaces include all separable Hilbert spaces,  all real hyperbolic spaces \cite{lyons2014}, 
and all open hemispheres \cite{lyons2019}. On some related information see \cite{b2016}. 

In  \cite{srb07} an alternative formula for ${\dCov}^2(X,Y)$ was given
in terms of characteristic functions $f_{X,Y}$, $f_X$ and $f_Y$ of $(X,Y)$, $X$, 
and $Y$ respectively. If the random variable $X$ takes values in a $p$-dimensional 
Euclidean space $\mathbb R^p$ and $Y$ takes values in $\mathbb R^q$ and both
variables have finite expectations we have 
\[
{\dCov}^2(X,Y):= \frac{1}{c_pc_q} \int_{\mathbb{R}^{p+q}} \frac{| f_{X,Y}(t,s)-
  f_X(t)f_Y(s)|^2}{|t|_p^{1+p}\;|s|_q^{1+q}} \,dt\,ds.  
\]
where $c_p$ and $c_q$ are constants. This formula clearly shows that
independence of $X$ and $Y$ is equivalent to $\dCov(X,Y) = 0$. 
On a generalization to dependence measures for more than two random vectors see \cite{bks2019}. 

If the expectations of $X, Y$ do not exist, we can generalize 
distance correlation for random variables with finite moments of order $\alpha >0$,
see \cite{srb07,sr09}. It is easy to see that the population
 distance correlation, $\dCor(X,Y)$, satisfies axioms (\ref{nii})
and (\ref{niv}). For the proof that  $\dCor(X,Y)$ satisfies
(\ref{ni}) and (\ref{niii}), see \cite{srb07}.  

An important generalization of distance correlation is \cite{S13}. This is 
related to a generalized distance correlation where the distance 
is a more general metric than the Euclidean one. These generalizations 
under some natural conditions like scale invariance also satisfy our axioms. 

In \cite{ms2019} we proved the following theorem which shows that 
in our axioms similarity cannot be replaced by stronger 
invariances like affine invariance (except in case $\dim H = 1$).

\begin{theorem} \label{th1}
Suppose $S$ is a set of pairs of nondegenerate random variables, and if 
$(X,Y) \in S$ then $\big(f(X), g(Y)\big) \in S $ for all affine transformations of $H$.  

If the dependence measure $\Delta(X,Y)$ on $S$ is invariant with
respect to all affine transformations $f,\,g$ of $H$ where $\dim H>1$ then
axiom (\ref{niv}) cannot hold. If $\dim H = 1$ then affinity is the  same as 
similarity and in this case distance correlation is affine 
invariant.  On the other hand, if $\Delta(X,Y)$ is invariant with
respect to all \textup{1--1} Borel measurable functions of $H$ then 
even if $\dim H = 1$, axiom (\ref{niv}) cannot hold.  
\end{theorem}

For an ``almost affine invariant'' version of distance correlation see \cite{degr2014}.
If we want to generalize the axioms from Hilbert spaces $H$ to general metric spaces 
$(\mathcal M, \delta)$ then  first we need a general definition of similarity in metric spaces.

\begin{definition}\label{defsim}
A mapping $f: \mathcal M \to \mathcal M$ is a {\it similarity} if 
there exists a constant $c > 0$ such that for all $x \in \mathcal M, y \in \mathcal M $ we have 
$\delta\big(f(x), f(y)\big) = c\delta(x,y) $.  
\end{definition}

A reformulation of our axioms to arbitrary metric spaces is the following.

\begin{enumerate}[(a)]
\item\label{dni} $\Delta (X, Y) = 0$ if and only if $X$ and $Y$ are
  independent. 
\item\label{dnii} $\Delta (X, Y)$ is invariant with respect to all
  similarity transformations of $(\mathcal M, \delta)$; that is, 
  $\Delta\big(f(X), g(Y)\big) =  \Delta(X,Y)$ where $f$, $g$ are similarity transformations of 
  $(\mathcal M, \delta)$. 
\item\label{dniii} $\Delta (X, Y) = 1$ if and only if $Y = f(X)$  with
  probability $1$, where $f$ is a similarity transformation of $(\mathcal M, \delta )$.
\item\label{dniv} 
  $\Delta (X, Y)$ is continuous; that is, if for some
  positive constant $K$ and $x_0\in \mathcal M, y_0 \in \mathcal M$ we have 
  $\E\big(\delta^2(X_n, x_0)+\delta^2(Y_n, y_0)\big)\le K$, $n=1,2,\dots$
  and $(X_n, Y_n)$ converges weakly (i.e., converges in distribution) to $(X,Y)$ 
  then $\Delta (X_n, Y_n) \to \Delta(X,Y)$.
  \end{enumerate}
  
Again, the condition on the boundedness of second 
  moments can be replaced by any other condition that guarantees the convergence of expectations: 
  $\E\delta(X_n , x_0) \to \E\delta(X, x_0) $ and $\E\delta(Y_n, y_0) \to \E\delta(Y, y_0)$; 
  the uniform integrability of $\delta(X_n,x_0), \delta(Y_n,y_0)$, which 
  follows from the boundedness of second moments, would equally do.

We will also need the following weaker forms of axioms (\ref{dnii}) and (\ref{dniii}):

\begin{enumerate}[(a*)]
\setcounter{enumi}{1}
\item\label{bb} $\Delta(X, Y)=\Delta(f(X), f(Y))$ for every similarity transformation $f$  
of $(\mathcal M, \delta)$.
\item\label{cc} $\Delta (X,Y) = 1$ if $Y=f(X)$ with probability $1$, where $f$ is a 
similarity transformation of $(\mathcal M, \delta)$.
\end{enumerate}

Our goal is to find a dependence measure that satisfies these axioms in arbitrary metric spaces. 

\section{How far can we go with distance correlation? }
\label{DC}

Distance correlation can be generalized to metric spaces $(\mathcal M, \delta)$ 
that are of negative type \cite{lyons2013}. A  metric space $(\mathcal M, \delta)$ is
called of negative type if the metric possesses the ``conditional negative definite'' property, namely
that  for all integers $n \ge 1$ and for all sets of $n$  points $x_i \in \mathcal M $ and 
$x'_i \in \mathcal M $ ($i = 1,2,\dots, n)$ and for all real numbers $a_1, a_2, \dots, a_n$ 
such that their sum is 0 we have
\[
\sum_{i,j} a_i a_j \delta (x_i, x_i') \le 0.
\]
Strong negative type metric spaces satisfy this with equality iff $a_1 = \dots =a_n = 0$. 
However, for the strong negative type property we need somewhat more, namely for all probability 
measures $\mu$ and $\nu$ defined on the Borel sets of $\mathcal M$
\[
\int \delta(x,y)d(\mu - \nu)^2(x,y) \le 0
\]
with equality iff $\mu = \nu $. 

According to a classical theorem of Schoenberg \cite{Sch37, Sch38} a necessary and 
sufficient condition for negative type of $(\mathcal M, \delta)$ is that 
$(\mathcal M, \sqrt{\delta})$ is isometrically embeddable into a Hilbert space. 
Obviously this property does not hold for every metric space. When it does then in these  
``nice'' metric spaces we can apply distance correlation, for all others we need to make new efforts. 

We can try to work with functions of $\delta$, say $\delta^*(\delta)$, that satisfies our 
axioms.  If the only problem is that the metric is not of strong negative type, only of 
negative type then it is easy to find a remedy: take the square root (or any other power 
$ 0 < r < 1$) of the metric and this new metric becomes of strong negative type, see 
\cite{lyons2013}.  

For arbitrary finite $\mathcal M$ we can show, see \cite{sr14}, that for a suitably large 
number $K$ the new distance $\delta^*(x,y) = \delta(x,y) + K$ whenever $x \neq y$ and $0$ otherwise, 
is always conditionally negative definite. On top of that, this simple transformation of the 
metric does not change the unbiased estimator of $\dCov$ which is simply invariant with 
respect to this additive constant $K$. 

For infinite $\mathcal M$  there does not always exist a strictly monotone increasing function 
$\delta^*(\delta) $ such that $(\mathcal M, \delta^*)$ is of negative type. 
Take e.g. two disjoint infinite sets, $A$ and $B$, and let $\mathcal M$ be their union.  Define the 
distance of two distinct elements to be $1$ if they are in different sets, and $2$ if
 they are in the same set. The function $\delta^*$  must have the following form: 
$\delta^*(1) = u $ , $\delta^*(2) = v$, $0 < u < v$. Define $a_i :=1$ for $n$ elements 
of $A$ and $a_i := -1$ for $n$ elements of $B$. Then the sum we need to check is  
$n(n-1)v - n^2 u$, which is positive for large enough $n$.

Another approach is this. If all we want from our dependence measure is to test independence
 then it is acceptable to change the distances in $(\mathcal M, \delta)$ and thus change the 
 distance correlation so long as we do not change $\dCor(X,Y) = 0$. If $f$ is an arbitrary \textup{1--1} 
Borel function on $(\mathcal M, \delta)$ and $X,Y$ are $(\mathcal M, \delta)$ valued 
random variables then they are independent iff $f(X), f(Y)$ are independent.  But every metric 
space is Borel isomorphic to a ``nice'' metric space that is embeddable isomorphically 
into a Hilbert space. According to Kuratowski's theorem two complete separable Borel spaces 
are Borel isomorphic iff they have the same cardinality. They are Borel isomorphic either to 
$\mathbb R$, or to $\mathbb Z$ or to a finite metric space. Denote this Borel isomorphism 
by $f$. If we can construct it then we can check the independence of the real valued random 
variables $f(X), f(Y)$ via distance correlation and this is equivalent 
to testing the independence of $X,Y$ that take values in general metric spaces. 
We might want to make $f$ continuous to avoid the negative effect of minor noise. In this case 
we can choose $f$ to be a homeomorphism between our metric space and a subspace of a 
Hilbert cube. This $f$ exists if and only if our metric space is separable. Here is how to construct 
such an $f$. 
 
Assume $\delta \le 1$ (otherwise, use $\delta/(\delta+1)$). Choose a dense countable 
sequence $(x_n)$ from $\mathcal{M}$ which exits because the metric space is separable, and define $f(x) := 
(\delta (x, x_n)/n)_{n\ge 1}$, a point in the Hilbert cube and here we can apply distance correlation 
for testing independence. 

These tricks can help to solve some of the problems in testing independence but they do not solve 
the problem of finding a general measure of dependence applicable to general metric space valued 
random variables. To use one of John von Neumann's favorite expressions, our goal here is to define
a dependence measure that applies to the ``rest of the universe''. 

\section{The population value of the earth mover's correlation}
\label{W}

First of all recall the definition of the earth mover's distance for probability measures 
$\mu, \nu$ on general metric spaces $(\mathcal M, \delta )$.  
We suppose that the topology of this metric space and the probability measures on
 the  Borel sets are ``compatible'', that is, we suppose that the probability measures are 
 Radon measures (finite on compact sets, outer regular and inner regular). 

Heuristically, if we have two (Radon) probability distributions, $\mu$ and $\nu$ on 
$(\mathcal M, \delta )$ then the earth mover's distance is the minimum cost of 
turning one pile of dust or dirt with distribution $\mu$ into the other with distribution $\nu$. 
The cost is proportional to the transport distance and also to the amount of dirt we transport. 

This distance was considered by \cite{m1781}, \cite{k1942},  \cite{kr1958} ,  
 \cite{wass1969}, \cite{rr1998}, \cite{v2009}, and many others, and in mathematical circles it is typically 
 called Wasserstein distance. Most statisticians and computer scientists call it earth mover's distance. 
 On a recent survey see \cite{PZ}.  On some recent advances see \cite{abm}, \cite {rtg}, and \cite {srgb}. 

Denote by $\mathcal P(\mathcal M)$ the set of all (Radon) probability measures $\mu$ on $\mathcal M$.
Suppose that for some $x_0 \in \mathcal M$ we have
\[
\int_{\mathcal M} \delta(x, x_0) d \mu(x) < + \infty.
\]
Then the earth mover's distance or Wasserstein distance of the probability measures $\mu$ and $\nu$ 
can be equivalently defined as
\[
e(\mu, \nu) : = \inf_{\gamma\in\Gamma(\mu,\nu)} \int_{\mathcal M \times \mathcal M} \delta (x, y)
 d\gamma (x,y), 
\]
where $\Gamma(\mu,\nu)$ is the set of all possible couplings of probability measures $\mu$ and $\nu$,
that is, the set of all joint distributions $\gamma$ of $(X,Y)$ with marginal distributions $\mu$ and $\nu$,
respectively. Equivalently, 
\[
e(\mu, \nu) =  e(X,Y):= \inf_{\gamma\in\Gamma(\mu,\nu)} \E[\delta(X,Y)],
\]
where again the infimum is taken for all joint distributions of $(X,Y)$ with marginal 
distributions $\mu$ and $\nu$, respectively.  

Mathematically this is not an easy minimization problem to solve. Even if $(\mathcal M, \delta )$
is an Euclidean space where the transportation cost is the Euclidean distance the solution is related
to the so-called Monge--Amp\`ere difference equation \cite{bb2000,c2003,cm2010}. 
For real valued random variables $X$, $Y$, however,  there is a simple formula for the earth mover distance. Denote $F(x)=\Pr(X\le x)$ and 
$G(y)=\Pr(Y\le y)$  the cdf's of $X$ and $Y$ and consider their generalized inverses $F^{-1}(u)$, $G^{-1}(u)$, 
defined as $F^{-1}(u)=\sup\{t: F(t)\le u\}$. Then
\begin{multline*}
e(X,Y)=\E\big|F^{-1}(U)-G^{-1}(U)\big|\\
=\int_0^1\big|F^{-1}(u)-G^{-1}(u)\big|du=\int_{-\infty}^{\infty}\big|F(t)-G(t)\big|dt.
\end{multline*}

Define a metric $d$ on the space $\mathcal M \times \mathcal M$, e.g. $d$ can be the Manhattan
 distance: $d\big[(x,y), (u,v)\big] = \delta(x,u) + \delta (y,v) $. 

\begin{definition}\label{emcov}
The earth mover's covariance of random variables $X,Y$ taking values in $(\mathcal M, \delta )$ 
is  the earth mover's distance between the joint distribution and the product of its 
marginals:
\[
\eCov(X,Y) = \inf_{\gamma\in\Gamma} \E d\big[(X,Y), (X', Y')\big] = e\big[(X,Y),  (X', Y')\big],
\]
where $\Gamma$ is the set of all possible joint distributions of the random variables $X,Y,X',Y'$  such
that $X'$ and $X$ are identically distributed, $Y'$ and $Y$ are also identically distributed, and 
$X', Y'$ are independent (and the joint distribution of $X$ and $Y$ is given).
\end{definition}

In the following we do not really need that $d$ is a Manhattan distance, what we need is more 
general, namely that $(\mathcal M \times \mathcal M , d)$ with a metric $d$ is a metric space 
such that 
\[
d\big[(x,u),(x,v)\big]=\delta(u,v),\ 
d\big[(x,u),(y,u)\big]=\delta(x,y),\ 
d\big[(x,x),(u,v)\big]\ge\delta(u,v).
\]

The following inequality is of Cauchy--Bunyakovsky--Schwarz type.
\[
e^2\big[(X,Y),(X',Y')\big]\le e\big[(X,X),(X,X')\big]\,e\big[(Y,Y),(Y,Y')\big],
\]
where $X$ and $X'$ are iid, as well as $Y$ and $Y'$, and $X'$, $Y'$ are independent.

In fact, we can show more, namely that 

\begin{theorem}\label{prop1}
\begin{equation}\label{ineq}
e\big[(X,Y),(X',Y')\big]\le \min\left\{e\big[(X,X),(X,X')\big],\,
e\big[(Y,Y),(Y,Y')\big]\right\}.
\end{equation}
\end{theorem}

\begin{proof} 
Suppose that the 
right-hand side is equal to $e\big[(Y,Y),(Y,Y')\big]$. 
In the sequel all random 
variables denoted by $X$ with or without subscripts or superscripts will be equidistributed 
with $X$, and the same holds for $Y$. Let $Y_2$ and $Y_3$ be 
independent, then 
\[
\E d\big[(Y_1,Y_1),(Y_2,Y_3)\big]\ge \E\delta(Y_2,Y_3)= \E d\big[(X_2,Y_2),(X_2,Y_3)\big],
\]
where $X_2$ is chosen in such a way that $(X_2,Y_2)$ and $(X,Y)$ are identically distributed, 
and $Y_3$ is independent of $(X_2,Y_2)$. Then the right-hand side is greater than or equal to 
$e\big[(X,Y),(X',Y')\big]$, while  the infimum of the left-hand side as $Y_1$ varies is just 
$e\big[(Y,Y),\allowbreak (Y,Y')\big]$.\qed
\end{proof}

On the right-hand side of \eqref{ineq} $e\big[(X,X),(X,X')\big]=\eCov(X,X)$ will be called \emph{the earth 
mover's variance}.

\begin{definition}\label{emvar}
The earth mover's variance of the distribution of $X$ is
\begin{equation}\label{evar}
\eVar(X):= \eCov(X,X) = e\big[(X,X),(X,X')\big].
\end{equation}
\end{definition}

\begin{theorem}\label{eVar}
The earth mover variance is the same as Gini's mean difference:
\begin{equation}\label{Gini}
\eVar(Y)=\E\delta(Y,Y'), 
\end{equation}
where $Y$ and $Y'$ are iid. 
\end{theorem}
\begin{proof}
We have seen above that
\[
\eVar(Y)=\inf_{Y_1}\E d\big[(Y_1,Y_1),(Y_2,Y_3)\big]\ge \E\delta(Y_2,Y_3),
\]
and equality is attained for $Y_1=Y_2$.
\end{proof}

\begin{example}  Let $X$ be an iid sample of size $n$ from the 
uniform distribution $U[0;1]$, apply the Euclidean metric in $\mathbb R^n$ and the 
Manhattan distance for pairs. Then by Remark \ref{rem1} below $\eVar(X)=\E|X-X'|$, where $X$ and $X'$ are 
independent uniform random points of the $n$ dimensional unit cube. For $n=1$ 
we get $\eVar(X)=1/3$. For general $n$ it is known that using the notation 
$\operatorname{\mathrm{erf}}(u)$ for the  ``error function'', i.e. the integral from $-u$ 
to $u$ of the Gaussian probability density function with $0$ expectation and variance $1/2$ we have
\[
\eVar(X)=\frac{1}{\sqrt{\pi}} \int_{0}^{\infty} \left \{1 - \left( \frac{\sqrt{\pi} 
\operatorname{\mathrm{erf}}(u)}{u} - \frac{1 - e^{-u^2}}{u^2} \right)^{\!\!n} \right\} \frac{du}{u^2}.
\]
We do not know any simple analytic expressions for $\eVar(X)$ if $n$ is 
arbitrary. However, by the inequality $|X-X'|\ge n^{-1/2}\sum_{i=1}^n|X_i-X'_i|$ it easily 
follows that $\eVar(X)\ge{\sqrt n}/3$. On the other hand, since the diameter of the unit 
cube is $\sqrt{n}$, we clearly have $\eVar(X)\le\sqrt{n}$. 
A somewhat better upper estimate is
\begin{multline*}
\E|X-X'|\le\left[\E|X-X'|^2\right]^{1/2}=
\left[n \E(X_1-X'_1)^2\right]^{1/2}\\
=\left[2n\Var(X_1)\right]^{1/2}=\sqrt{n/6}.
\end{multline*}
\end{example}
Based on Theorem \ref{prop1} we can now introduce the definition of a new type of correlation. 

\begin{definition}\label{emcorr}
The \emph{earth mover's correlation} 
of the distributions of $X$ and $Y$ is defined  as
\[
\eCor(X,Y)=\frac{\eCov(X,Y)}{\min\big\{\eVar(X),\eVar(Y)\big\}}\,.
\]
We do not define $\eCor(X,Y)$ when $\min\big\{\eVar(X),\eVar(Y)\big\} = 0$.
\end{definition}
\begin{remark}\label{rem1} 
By the previous theorem in the formula for $\eCor$ the denominator $\min\big\{\eVar(X),\eVar(Y)\big\} 
= \min\big\{\E\delta(X,X'), \E\delta(Y,Y') \big\} = 0$ 
iff at least one of $X,Y$ is constant with probability $1$. In this case we do not define $\eCor$. 
It is interesting to note that for real valued random variables $\eVar$ is easy to compute.
It is known, see e.g. \cite{Yitzhaki}, that

\[
\eVar(X)=2\int_{-\infty}^{\infty}F(x)(1-F(x))\,dx,
\]
where $F(x)=\Pr(X\le x)$ is the cdf of the random variable $X$. 
\end{remark} 

\begin{remark}\label{rem2}
Let us apply the Manhattan distance for pairs. Then by the triangle inequality for $\delta$ we have
$\delta(X,X')+\delta(Y,Y')\ge |\delta(X,Y)-\delta(X',Y')|$, thus
\[
\eCov(X,Y)\ge \inf_{(X',Y')} \E\big|\delta(X,Y)-\delta(X',Y')\big|\ge\big|\E\delta(X,Y)-\E\delta(X',Y')\big|.
\]
\end{remark}

\begin{example}\label{ex2}
Let $X$ and $Y$ be indicators, $\Pr(X=1)=1-\Pr(X=0)=p_X$, $\Pr(Y=1)=1-\Pr(Y=0)=p_Y$, 
$\Pr(X=Y=1)=p_{XY}$. Let us apply the Euclidean metric in $\mathbb R$ and the Manhattan 
distance for pairs. Then
\[
\eCor(X,Y)=\frac{|p_{XY}-p_X p_Y|}{\min\big\{p_X(1-p_X), \,p_Y(1-p_Y)\big\}}\,.
\]

Indeed, in the lower bound of Remark \ref{rem2} we have
\begin{align*}
\E|X-Y|&-\E|X'-Y'|\\
&=\Pr(X\ne Y)-\Pr(X'\ne Y')\\
&=\Pr(X'=Y')-\Pr(X=Y)\\
&=p_X p_Y+(1-p_X)(1-p_Y)-p_{XY}-(1-p_X-p_Y+p_{XY})\\
&=2(p_X p_Y-p_{XY}),
\end{align*}
thus $\eCov(X,Y)\ge 2|p_{XY}-p_X p_Y|$. 

On the other hand, we will
construct random variables $X,Y,X',Y$ with the desired distribution in such a way that
\[
\E\big(|X-X'|+|Y-Y'|\big)=2|p_{XY}-p_X p_Y|.
\]

Let $U,V$ be independent and uniformly distributed on $[0,1]$, and define
\begin{gather*}
X=X'=I(U\le p_X),\quad Y'=I(V\le p_Y),\\
Y=I\left(U\le p_X,\  V\le\frac{p_{XY}}{p_X}\right)+I\left(U>p_X,\ V\le\frac{p_Y-p_{XY}}{1-p_X}\right)
\end{gather*}
Then $\Pr(X\ne X')=0$ and $\Pr(Y\ne Y')=\Pr(Y=1,Y'=0)+\Pr(Y=0,Y'=1)$. Here
\begin{multline*}
\Pr(Y=1,Y'=0)\\=\Pr\left(U\le p_X,\ p_Y<V\le \frac{p_{XY}}{p_X}\right)+
\Pr\left(U>p_X,\ p_Y<V\le\frac{p_Y-p_{XY}}{1-p_X}\right),
\end{multline*}
and similarly,
\begin{multline*}
\Pr(Y=0,Y'=1)\\=\Pr\left(U\le p_X,\ \frac{p_{XY}}{p_X}<V\le p_Y\right)+
\Pr\left(U>p_X,\ \frac{p_Y-p_{XY}}{1-p_X}< V\le p_Y\right).
\end{multline*}
Altogether we have
\[
\Pr(Y\ne Y')=p_X\Big|p_Y-\frac{p_{XY}}{p_X}\Big|+
(1-p_X)\Big|p_y-\frac{p_Y-p_{XY}}{1-p_X}\Big|=2|p_{XY}-p_X p_Y|,
\]
thus $\eCov(X,Y)=2|p_{XY}-p_X p_Y|$.

Finally, $\eVar(X)=2p_X(1-p_X)$ is straightforward, a special case of the previous formula.
\end{example}

The absolute value of Pearson's correlation $\rho$ for indicators is

\[
|\rho(X,Y)|=\frac{|p_{XY}-p_X p_Y|}{\sqrt{p_X(1-p_X)\strut \,p_Y(1-p_Y)}}\,
\]
thus for indicators $X$ and $Y$ we have $|\rho(X,Y)| \le  \eCor(X,Y)$ (and we have equality iff $p_X = p_Y$). 

Based on this observation one can suspect that $|\rho(X,Y)| \le  \eCor(X,Y)$ for all real valued random variables 
with finite variance. This conjecture is also supported by the fact that the independence of $X,Y$ implies 
their uncorrelatednes. 
In the other extreme case when $\rho(X,Y) = \pm 1$ we know that $Y = f(X)$ where 
$f$ is a similarity (here a linear function) and by Theorem \ref{th3} below in this case we have 
$\eCor(X,Y) = 1$. 

The conjecture that $|\rho(X,Y)| \le  \eCor(X,Y)$ holds for all real valued random variables 
with finite variance, however,  can easily be disproved. The following theorem shows that 
if the joint distribution of $X,Y$ is bivariate normal, the opposite inequality holds. 

\begin{theorem} \label{thn}
Let $(X,Y)$ be bivariate normal with correlation $\varrho(X,Y)=\varrho$. Then 
\[
\eCor(X,Y)\le\left[1-\sqrt{1-\varrho^2}\,\right]^{1/2}\le |\varrho|,
\]
and the last inequality is strict unless $\varrho=0$ or $\varrho= \pm 1$.
\end{theorem}
Actually we  have the following 

\begin{conjecture}\label{Conj1}
Let $(X,Y)$ be bivariate normal with correlation $\varrho(X,Y)=\varrho$. Then
$\eCor(X,Y) = \left[1-\sqrt{1-\varrho^2}\,\right]^{1/2}$.
\end{conjecture}

\noindent \textbf{Proof of Theorem \ref{thn}.}
Let $(X,Y)$ be bivariate normal with $\Var(X)=\sigma_X^2$, $\Var(Y)=\sigma_Y^2$ and
$\varrho(X,Y)=\varrho$. We can suppose $\E X=\E Y=0$ and $\sigma_X^2\ge\sigma_Y^2$.  
Let $X'$ and $Y'$ be independent zero mean 
normal with variances $\sigma_X^2$ and $\sigma_Y^2$, respectively. Finally, set $X=X'$ and 
$Y=(\sigma_Y/\sigma_X)\varrho X'+\sqrt{1-\varrho^2}\,Y'$. Then $X,Y$ have the prescribed
joint distribution, and $Y-Y'$ is normal with mean $0$ and variance 
\[
\sigma_Y^2\left[\varrho^2+\left(1-\sqrt{1-\varrho^2}\,\right)^{\!2}\right]=
2\sigma_Y^2\left[1-\sqrt{1-\varrho^2}\,\right],
\]
hence
\[
\eCov(X,Y)\le\E|Y-Y'|=\frac{2\sigma_Y}{\sqrt{\pi}}\left[1-\sqrt{1-\varrho^2}\,\right]^{1/2}.
\]
In the denominator of $\eCor$ we have $\eVar(Y)=\dfrac{2\sigma_Y}{\sqrt{\pi}}\le 
\dfrac{2\sigma_X}{\sqrt{\pi}}=\eVar(X)$,
thus
\[
\eCor(X,Y)\le\left[1-\sqrt{1-\varrho^2}\,\right]^{1/2}\le |\varrho|,
\]
and the last inequality is strict unless $\varrho=0$ or $\varrho= \pm 1$.\qed

Concerning the lower bound of $\eCor(X,Y)$, if $\sigma_X = \sigma_y$ then Remark \ref{rem2} 
provides the following inequality:
\[
\big|1 -\sqrt{1 - \varrho}\,\big| \le \eCor(X,Y).
\]

The arguments in the proofs support the next conjecture. 

\begin{conjecture}\label{Conj2}
In computing the infimum $\eCov(X,Y)=\inf_{(X',Y')}\E\big[\delta(X, X')+\delta(Y, Y')\big]$, under 
``general conditions'' we can suppose $X=X'$ or $Y=Y'$. 
\end{conjecture}

On the above mentioned ``general conditions''  see below. But first we show by an example that the conjecture is not true without some restrictions.

\begin{example}
If $X$ and $Y$  are \textup{1--1} functions of each other then the conjecture would imply that 
$\eCor(X,Y) = 1$ because $Y$ is a function of $X = X'$ thus $Y$ is independent of $Y'$.  Hence  
$\eCov(X,Y) = \min\{ \eVar(X), \eVar(Y)\}$. Thus in case of continuous marginals 
the empirical $\eCor$ would always be $1$ because for continuous marginals no vertical or horizontal lines 
can contain more than one sample points with probability one. This is, however, not true as is shown by the following sample of four elements:
$(1,4),\, (2,2),\, (3,3),\, (4,1)$. 
Here $\eVar = 5/4$ for both coordinates but $\eCov = 1$.
\end{example}

\begin{theorem}\label{conj}
Conjecture \ref{Conj2} implies Conjecture \ref{Conj1}.
\end{theorem}

\begin {proof}
The infimum in the theorem can be computed by applying conditional quantile transformations.
Suppose $X=X'$. Let $F(x)$, $G(y)$ denote the cdf of $X$ and $Y$, resp., and
$G(y|x)=\Pr(Y\le y\mid X=x)$, the conditional cdf of $Y$. Then the infimum of 
$\E|Y-Y'|$ under the condition that  $Y=Y'$ in distribution, but $X',Y'$ are independent, equals
\[
\E\bigg(\int_{-\infty}^{\infty}\big|G(y|X)-G(y)\big|dy\bigg)=
\int_{-\infty}^{\infty}\E\big|G(y|X)-G(y)\big|\,dy.
\]
Note that $G(y)=\E G(y|X)$, thus the integrand on the right hand side is a kind of a mean 
absolute difference. An alternative formula for $\eCov$ is
\[
\eCov(X,Y)=\int_{-\infty}^{\infty}\int_{-\infty}^{\infty}\big|G(y|x)-G(y)\big|\,dF(x)\,dy.
\]

In the case of jointly normal $X,Y$   the conditional quantile transformation leads to the same
representation of $Y$ as a linear combination of $X'$ and $Y'$ that we used in the proof 
of Theorem \ref{thn}. Thus our Conjecture \ref{Conj1} would follow from Conjecture \ref{Conj2}.
\qed
\end{proof}
Unfortunately we could not find simple ``general conditions'' for the validity of Conjecture \ref{Conj2}.

It is easy to see that $\eCor$ as a new measure of dependence satisfies at least two
of our axioms for dependence measures. 
Axioms (\ref{dni}), and (\ref{dniv}) hold. Concerning (\ref{dnii}) and (\ref{dniii}) we can only prove 
the weaker (\ref{bb}*) and (\ref{cc}*).

\begin{theorem}\label{th3}\ %

$\eCor(X, Y)=\eCor(f(X), f(Y))$ for every similarity transformation $f$ of our metric space.

If $Y = f(X)$ where $f$ is a similarity transformation then $\eCor(X, Y) = 1$.
\end{theorem}

\begin{proof}
From the definition it is obvious that $\eCov(f(X),f(Y))=c\cdot\eCov(X,Y)$. Therefore we also have 
$\eVar(f(X)) = c\cdot \eVar(X)$, and finally $\eCor(f(X),f(Y))=\eCor(X,Y)$.
 
For independent $X_1, X_2, X_3$ we have 
\begin{align*}
d\big[(X_1,f(X_1)),(X_2,f(X_3))\big]&=\delta(X_1,X_2)+\delta(f(X_1),f(X_3))\\
&=\delta(X_1,X_2)+c\cdot\delta(X_1,X_3)\\
&\ge \min\{1,c\}\big[\delta(X_1,X_2)+\delta(X_1,X_3)\big]\\
&= \min\{1,c\}\, d\big[(X_1,X_1),(X_2,X_3)\big]\\
&\ge\min\{1,c\}\eVar(X)\\
&=\min\{\eVar(X),\eVar(f(X))\}.
\end{align*}
The infimum of the left hand side as $X_2$ and $X_3$ remain independent is equal to $\eCov(X, f(X))$. Thus 
$\eCor(X, f(X)) \ge 1$. The other direction follows from Theorem \ref{prop1}. \qed
\end{proof}

Thus we proved the following result.  

\begin{theorem}\label{th4}
 In arbitrary metric spaces $(\mathcal M , \rho)$ the earth mover's correlation  $\Delta (X, Y) = \eCor(X,Y)$ 
satisfies axioms  (\ref{dni}), (\ref{bb}*),  (\ref{cc}*), and (\ref{dniv}).\qed
\end{theorem}

It is easy to see that for an arbitrary metric space 
$(\mathcal M, \delta)$  it cannot be true that $\eCor(X,Y) = 1$ always implies  $Y= f(X)$ 
where $f$ is a similarity. A counterexample is the following. Let $\mathcal M$ be the set of points 
of the Euclidean plane with the usual Euclidean metric. Suppose that here  $\eCor(X,Y) = 1$ implies  $Y= f(X)$ 
where $f$ is a similarity. If the random variables $X$ and $Y$ are supported on the $x$ line then we know 
that the similarity is $Y = aX +b$. 
Now define a new metric on the plane as follows: $\delta(x,y) = |x - y|$ if both 
$x$ and $y$ are on the $x$ coordinate axis (the second coordinate is $0$), otherwise for all $x\neq y$ 
define $\delta(x,y) = |x-y| + 1$. This does not change $\eCor(X,Y) = 1$ because $X$ and $Y$ are 
supported on the $x$ but $y = ax +b$ cannot be extended to the whole plane as a similarity 
with respect to the new metric.

\begin{conjecture}
 For Banach space valued random variables we have the iff statement in axiom (\ref{dniii}): 
$\eCor (X, f(X)) = 1 $ if and only if $Y = f(X)$ with probability $1$, where $f$ is a similarity 
transformation of the Banach space. 
\end{conjecture}

Although we could not prove this conjecture it is interesting to note that by a theorem of 
\cite{mu1932}, any bijective similarity $f$ of any Banach space (or of any normed linear space) 
is affine, that is, $f(x) - f(0)$ is linear. Thus similarities in Banach spaces must have a very simple structure. 

By the way, it is interesting to note that we can always embed every metric space 
$(\mathcal M, \delta)$ into the Banach space $C_b(\mathcal M)$  of bounded continuous 
functions on $(\mathcal M, \delta)$, just take the function 
\[
f(x)(y) := \delta(x,y) - \delta(x_0, y),
\]
where $x_0$ is an arbitrary element of $\mathcal M$. 

We note that one can easily define the earth mover's correlation for more than
two variables. The population version of  $\eCov$ for three variables is as follows:
\[
\eCov(X,Y,Z) = \inf_{(X',Y',Z')} \E d\big[(X,Y,Z), (X', Y', Z')\big].
\]
Here in distribution $X= X'$, $Y=Y'$, $Z=Z'$, and $X',Y',Z'$ are independent, and we take the $\inf$ 
over all joint distributions of $(X,Y,Z)$ and $(X',Y',Z')$.

The population version of the three-variate earth mover's correlation is
\[
\eCor(X,Y, Z)=\frac{\eCov(X,Y,Z)}{\min\big\{\eVar(X), \eVar(Y), \eVar(Z)\big\}}\,.
\]

Thus we have a natural measure for mutual dependence of more than two random variables.

\section{Empirical earth mover's correlation }
\label{E}

The earth mover's metric suggests the following earth mover's distance definition between two 
sequences $x:= (x_1,x_2,\dots, x_n)$ and $y:=(y_1,y_2,\dots, y_n)$:
\[ 
\mathcal{E}(x,y) := \inf_{\pi} \sum_{i=1}^n \delta(x_i , y_{\pi(i)}),
\]
where the infimum is taken for all permutation $\pi$ on the integers $1,2,\dots, n$.
One can easily see that for real valued data, if the ordered sample is denoted by subscripts in 
brackets, then 
\[
\mathcal{E}(x,y) := \sum_{i=1}^n |x_{(i)} - y_{(i)}|.
\]

The empirical version of $\eCov$ is the minimum transportation cost between the following two 
mass distributions or probability distributions:

($Q_1$)   $1/n$ mass at each point  $(x_i, y_i) ,\  i = 1,2,\dots, n$ 

\noindent and

($Q_2$) $1/n^2$  mass at each point $(x_i, y_j) ,\  i,j = 1,2,\dots, n$.

It is easy to see that the empirical $\eVar$ is the arithmetic average of the distances $\delta(x_i , x_j)$ 
because the cost to transport $1/n^2$ mass from the point $(x_i, x_j)$ to the main diagonal $(x,x)$
is at least $\delta (x_i , x_j)/n^2$ and we can achieve this via ``horizontal''  transportation only. 
This is not the case if we want to transport to $n$ general points, not necessarily 
on the main diagonal. The ``naive'' computational complexity of the empirical $\eVar$ which is 
essentially Gini's mean difference is $O(n^2)$ but for real valued random variables we can decrease 
it to $O(n \log n)$.

The complexity of the computation of the empirical $\eCov$ is less obvious. 

Our transportation problem can be reduced to an assignment problem between two sets of $n^2$ 
points thus according to the ``Hungarian algorithm'' \cite{kuhn1955} this optimization can be 
solved in polynomial time. It was shown by \cite{ek72} and \cite{tomi} that the algorithmic 
complexity of assignment problem for two sets of $n$ points is $O(n^3)$ 
thus in our case the complexity can be reduced to $O(n^6)$.
 
This is not very encouraging. A better complexity, namely $O(n^3 \log^2 n)$,  is in \cite{ks1995}. 
Here the authors show that for the (linear) transportation problem with $m$ supply nodes, $n$ demand 
nodes and $k$ feasible arcs there is an algorithm which runs in time proportional to  
$m \log m(k + n \log n)$  assuming w.l.o.g. that $m \ge n$,  still at least one order of magnitude 
worse than the algorithmic complexity, $O(n^2)$, of computing the distance covariance or the distance correlation. 
This is the price we need to pay for the generality of $\eCov$ and $\eCor$.  The AMPL 
(A Mathematical Programming Language) code is easy to apply for computing empirical $\eCov$ 
and then $\eCor$. In \cite{akt1984} it was shown that given $n$ random blue and $n$ random 
red points on the unit square, the transportation cost between them is typically $\sqrt{n \log n}$. 
Our problem is to find the optimal transportation costs when the distance is the Manhattan distance 
and the number of red points is different from the number of blue points (the total mass is the 
same). A recent paper \cite{afpvx2017} suggests that our task of computing the earth mover's
distance between two sets of size $n^2$ can be done with the first algorithm in the cited paper with
 $O(\log^2(1/\varepsilon))$ approximation error bound in $O(n^{2 +\varepsilon})$ steps, for 
 any $\varepsilon > 0$. On related algorithmic optimizations see \cite{AWR} and \cite{ABRW}. 
 
 \section{Conclusion}
 For Hilbert space valued random variables in \cite{ms2019} we proved that distance correlation
 is a good mesure of dependence in the sense that distance correlation satisfies our axioms 
 (\ref{ni})--(\ref{niv}). For general metric space valued random variables, however, this is not true. 
 The earth mover's correlation ($\eCor$) introduced in this paper works for general metric spaces 
 in the sense that $\eCor$ satisfies axioms (\ref{dni}), (\ref{bb}*), (\ref{cc}*), (\ref{dniv}), and we 
 conjecture that under general conditions, e.g. for Banach space valued random variables, $\eCor$ 
 satisfies (\ref{dni}), (\ref{dnii}), (\ref{dniii}), (\ref{dniv}), too. These are counterparts of axioms 
 (\ref{ni})--(\ref{niv}). Our main result is Theorem \ref{prop1}, the earth mover's version of the 
 Cauchy--Bunyakovsky--Schwarz inequality. Conjectures \ref{Conj1} and \ref{Conj2} are challenges 
 for further research aiming easier computations of $\eCor$. If all we want is to test independence 
 then we do not really need the empirical $\eCor$, it is simpler to work with the empirical earth mover's
 covariance. For similar statistical tests see \cite{DS, GS}.

\vspace{2cm}

\noindent\textbf{Tam\'as F. M\'ori}\\
Alfr\'ed R\'enyi Institute of Mathematics\\
Re\'altanoda u. 13--15.\\
H-1053 Budapest, Hungary\\
{\tt mori.tamas@renyi.hu}\\

\noindent\textbf{G\'abor J. Sz\'ekely}\\
National Science Foundation\\
2415 Eisenhower Avenue\\
Alexandria, VA 22314, USA\\
{\tt gszekely@nsf.gov}


\vspace{-3ex}
\begin{thebibliography}{99}
\setlength{\parskip}{-3pt}\vspace{-2ex}

\bibitem{afpvx2017}
\textbf{Agarwal, P. K., K. Fox, D. Panigrahi, K. R. Varadarajan and A. Xiao,}
Faster algorithms for the geometric transportation problem,
in: B. Aronov and M. J. Katz (Eds.) \textit{33rd International Symposium on Computational 
Geometry (SoCG 2017)}, Leibniz International Proceedings in Informatics (LIPIcs) \textbf{77}, 
Schloss Dagstuhl--Leibniz-Zentrum f\"ur Informatik, Dagstuhl, Germany, 2017, 7:1--7:16.

\bibitem{akt1984}
\textbf{Ajtai, M., J. Koml\'os J and G. Tusn\'ady,}
On optimal matchings,
\textit{Combinatorica}, \textbf{4(4)} (1984), 259--264. 

\bibitem{AWR}
\textbf{Altschuler, J., J. Niles-Weed and P. Rigollet,}
Near-linear time approximation algorithms for optimal transport via Sinkhorn iteration. 
In: I. Guyon et al. (Eds.) \textit{Advances in Neural Information Processing Systems 30 (NIPS 2017)},
Curran Associates, Inc., Red Hook, NY, 2017, 1964--1974.

\bibitem{ABRW}
\textbf{Altschuler, J., F. Bach, A. Rudi and J. Weed,}
Massively scalable Sinkhorn distances via the Nystr\"om method.
arXiv:1812.05189v3 [stat.ML] preprint,
\url{https://arxiv.org/pdf/1812.05189}, Accessed 21 January 2020. (2019)

\bibitem{abm}
\textbf{Anderes, E., S. Borgwardt and J. Miller,}
Discrete Wasserstein barycenters: optimal transport for discrete data,
\textit{Math. Meth. Oper. Res.}, \textbf{84(3)} (2016), 389--409.

\bibitem{bb2000}
\textbf{Benamou, J. D.  and Y. Breiner,}
A computational fluid mechanics solution to the Monge--Kantorovich mass transfer problem,
\textit{Numer. Math.}, \textbf{84(3)} (2000), 375--393.

\bibitem{b2016}
\textbf{Bingham, N. H., A. Mijatovi\'c and T. L. Symons,} 
Brownian manifolds, negative type and geo-temporal covariances,
\textit{Commun. Stoch. Anal.}, \textbf{10(4)} (2016), 421--432. 

\bibitem{bks2019}
\textbf{B\"ottcher, B., M. Keller-Ressel and R. L. Schilling,}
Distance multivariance: New dependence measures for random vectors,
\textit{Ann. Statist.}, \textbf{47(5)} (2019), 2757--2789.

\bibitem{c2003}
\textbf{Caffarelli, L. A.,}
The Monge-Amp\`ere equation and optimal transportation, an elementary review, 
in: Ambrosio, L. et al. (Eds.) \textit{Optimal Transportation and Applications}, 
Lecture Notes in Mathematics \textbf{1813}, 
Springer, Berlin, Heidelberg, 2003, 1--10.

\bibitem{cm2010}
\textbf{Caffarelli, L. A. and R. J. McCann,}
Free boundaries in optimal transport and Monge-Amp\`ere obstacle problems,
\textit{Ann. of Math (2)}, \textbf{171(2)} (2010), 673--730.

\bibitem{DS}
\textbf{Deb, N. and B. Sen,}
Multivariate rank-based distribution-free nonparametric testing using measure transportation,
ArXiv 1909.08733 [math.ST] preprint.
\url{https://arxiv.org/pdf/1909.08733}, Accessed 06 June 2020 (2019)

\bibitem{degr2014} 
\textbf{Dueck, J., D. Edelmann, T. Gneiting and D. Richards,}
The affinely invariant distance correlation,
\textit{Bernoulli}, \textbf{20} (2014), 2305--2330.

\bibitem{ek72}
\textbf{Edmonds, J. and R. M. Karp,}
Theoretical improvements in algorithmic efficiency for network flow problems,
\textit{J. ACM}, \textbf{19} (1972), 248--264.

\bibitem{GS}
\textbf{Ghosal, P. and B. Sen,}
Multivariate ranks and quantiles using optimal transportation and applications to goodness-of-fit
testing,
ArXiv 1905.05340v2 [math.ST] preprint.
\url{https://arxiv.org/pdf/1905.05340v2}, Accessed 06 June 2020 (2019)

\bibitem{Jakobsen}
\textbf{Jakobsen, M. E.,} 
Distance Covariance in Metric Spaces: Non-Parametric Independence
Testing in Metric Spaces  (Master's thesis), arXiv:1706.03490 [math.ST] preprint,
\url{https://arxiv.org/pdf/1706.03490}, Accessed 21 January 2020 (2017)

\bibitem{k1942}
\textbf{Kantorovich, L. V.,} 
On the translocation of masses, 
\textit{Dokl. Akad. Nauk SSSR}, \textbf{37(7--8)} (1942), 227--229 (in Russian).

\bibitem{kr1958}
\textbf{Kantorovich, L. V. and G. S. Rubinstein,}
On a space of completely additive functions, 
\textit{Vestnik Leningrad Univ. Ser. Mat. Mekh. Astron.}, \textbf{13(7)} (1958), 52--59 (in Russian).

\bibitem{ks1995}
\textbf{Kleinschmidt, P. and H. Schannath,}
A strongly polynomial algorithm for the transportation problem,
\textit{Math. Program.}, \textbf{68} (1995), 1--13.

\bibitem{kuhn1955}
\textbf{Kuhn, H. W.,} 
The Hungarian method for the assignment problem,
\textit{Naval Res. Logist.}, \textbf{2} (1955), 83--97.

\bibitem{lyons2013}
\textbf{Lyons, R.,} 
Distance covariance in metric spaces,
\textit{Ann. Probab.}, \textbf{41(5)} (2013), 3284--3305.

\bibitem{lyons2014}
\textbf{Lyons, R.,} 
Hyperbolic space has strong negative type,
\textit{Illinois J. Math.}, \textbf{58(4)} (2014), 1009--1013.

\bibitem{lyons2018}
\textbf{Lyons, R.,} 
Errata to ``Distance covariance in metric spaces'',
\textit{Ann. Probab.}, \textbf{46(4)} (2018), 2400--2405.

\bibitem{lyons2019}
\textbf{Lyons, R.,} 
Strong negative type in spheres. ArXiv 1905.02863 [math.MG] preprint.
\url{https://arxiv.org/pdf/1905.02863}, Accessed 21 January 2020 (2019)

\bibitem{mu1932}
\textbf{Mazur, S. and S. Ulam,}
Sur les transformationes isom\'etriques d'espaces vectoriels norm\'es,
\textit{C. R. Acad. Sci. Paris}, \textbf{194} (1932), 946--948.

\bibitem{m1781}
\textbf{Monge, G.,} 
\textit{M\'emoire sur la th\'eorie des d\'eblais et des remblais},
De l'Imprimerie Royale, Paris, 1781.

\bibitem{ms2019}
\textbf{M\'ori, T. F. and G. J. Sz\'ekely,}
Four simple axioms of dependence measures,
\textit{Metrika}, \textbf{82} (2019), 1--16.

\bibitem{PZ}
\textbf{Panaretos, V. M. and Y. Zemel,}
Statistical aspects of Wasserstein distances,
\textit{Annu. Rev. Stat. Appl.}, \textbf{6(1)} (2019), 405--431.

\bibitem{rr1998}
\textbf{Rachev, S. T. and L. R\"uschendorf,}
\textit{Mass Transportation Problems: Volume I: Theory},
Springer-Verlag, New York, 1998.

\bibitem{rtg}
\textbf{Rubner, J., J. Tomasi and L. J. Guibas,}
The earth mover’s distance as a metric for image retrieval,
\textit{Int. J. Comput. Vis.}, \textbf{40(2)} (2000), 99--121.

\bibitem{Sch37}
\textbf{Schoenberg, I. J.,}  
On certain metric spaces arising from Euclidean spaces by a change of metric 
and their imbedding in Hilbert space,
\textit{Ann. of Math. (2)}, \textbf{38(4)} (1937), 787--793. 

\bibitem{Sch38}
\textbf{Schoenberg, I. J.,} 
Metric spaces and positive definite functions,
\textit{Trans. Amer. Math. Soc.}, \textbf{44(3)} (1938), 522--536. 

\bibitem{S13}
\textbf{Sejdinovic, D., B. Sriperumbudur, A. Gretton and K. Fukumiyu,}
Equivalence of distance-based and RKHS-based statistics in hypothesis 
testing,
\textit{Ann. Statist.}, \textbf{41} (2013), 2263--2291.

\bibitem{srgb}
\textbf{Solomon, J., R. Rustamov, L. Guibas and A. Butscher,}
Earth mover’s distances on discrete surfaces. 
\textit{ACM Trans. Graph.}, \textbf{33(4)} (2014), Article No.: 67

\bibitem{srb07}
\textbf{Sz\'ekely, G. J., M. L. Rizzo and N. K. Bakirov,}
Measuring and testing independence by correlation of distances,
\textit{Ann. Statist.}, \textbf{35} (2007), 2769--2794.

\bibitem{sr09}
\textbf{Sz\'ekely, G. J. and M. L. Rizzo,}
Brownian distance covariance,
\textit{Ann. Appl. Stat.}, \textbf{3} (2009), 1236--1265.

\bibitem{sr14}
\textbf{Sz\'ekely, G. J. and M. L. Rizzo,}
Partial distance correlation with methods for dissimilarities,
\textit{Ann. Statist.}, \textit{42} (2014), 2382--2412.

\bibitem{tomi}
\textbf{Tomizawa, N.,} 
On some techniques useful for solution of transportation network problems,
\textit{Networks}, \textbf{1} (1971), 173--194.

\bibitem{v2009}
\textbf{Villani, C.,} 
\textit{Optimal Transport: Old and New},
Springer-Verlag, Berlin, Heidelberg, 2009.

\bibitem{wass1969}
\textbf{Wasserstein, L. N.,} 
Markov processes over denumerable products of spaces describing large systems of automata,
\textit{Probl. Inf. Transm.}, \textbf{5(3)} (1969), 47--52.

\bibitem{Yitzhaki}
\textbf{Yitzhaki, S.,}  
Gini’s mean difference: a superior measure of variability for non-normal distributions,
\textit{Metron}, \textbf{61(2)} (2003), 285--316.

\end{thebibliography}
\end{document}